# The Donsker-Varadhan Action Functional for Probabilistic Cellular Automata

A. Eizenberg

*Azrieli College of Engineering, Jerusalem 9103501, Israel*
Email: alex@jce.ac.il

ABSTRACT

From the perspective of the large deviations theory of occupational measures, the paper considers Probabilistic Cellular Automata (PCA) as Markov chains on infinite dimensional space. It turns out that for a wide range of PCA, the corresponding Donsker-Varadhan action functional yields values other than infinity only on a narrow class of probability measures.

KEY WORDS: *Large Deviations; Probabilistic Cellular Automata; Markov Chains*

MATHEMATICAL SUBJECT CLASSIFICATION (2020) 60F10

## 1. INTRODUCTION

Let $X_t$, $t \in \mathbb{Z}^+$, be a time homogeneous Markov chain with a compact metric phase space $\Gamma$ equipped with the Borel $\sigma$-algebra $\mathcal{B}$. Denote by $M(\Gamma)$ the set of the probability Borel measures defined on $\Gamma$ equipped with the weak topology. Consider the sequence of the temporal occupational measures

$$\xi_T = \frac{1}{T}\sum_{t=0}^{T-1} \delta(X_t), \quad T \in \mathbb{Z}, T \geq 1, \tag{1.1}$$

where $\delta(x)$ is the unit measure concentrated at a point $x \in \Gamma$. Due to the classical results established by Donsker and Varadhan in [4] - [6], it is a well-known fact that under wide classes of conditions the asymptotic behavior of the occupational measures $\xi_T$ obeys the large deviations principle and can be described by means of the action functional $I: M(\Gamma) \to [0, \infty]$ defined for any $\nu \in M(\Gamma)$ by the formula

$$I(\nu) = -\inf\left\{\int_\Gamma \log\left(\frac{E_x f(X_1)}{f(x)}\right) \nu(dx): f \in \mathcal{U}\right\}, \tag{1.2}$$

where $\mathcal{U}$ is the set of positive continuous functions defined on $\Gamma$. In particular, if the Markov chain $X_t$ is a Feller process, then for any closed with respect to the weak topology subset $K$ of $M(\Gamma)$ the following upper bound holds uniformly with respect to $x \in \Gamma$,

$$\limsup_{T \to \infty} \frac{\log P_x(\xi_T \in K)}{T} \leq -\inf_{\nu \in K} I(\nu) \tag{1.3}$$

(see, for instance, [5] for the general case, or [8] for the case of Probabilistic Cellular Automata, where a simplified proof is provided). Moreover, it is known that under some additional conditions, such as, for instance, a Doeblin or Harris type conditions, or, more generally, certain uniformity conditions, such as, for example, Assumption U formulated in [7] , the following lower large deviations bounds have been derived:

$$\liminf_{T \to \infty} \frac{\log P_x(\xi_T \in U)}{T} \geq - \inf_{\nu \in U} I(\nu) \qquad (1.4)$$

for any open with respect to the weak topology subset $U$ of $M(\Gamma)$ and for each $x \in \Gamma$.

However, such assumptions, as a rule, are not satisfied for a large class of Markov chains, usually called Probabilistic Cellular Automata (PCA), which can be considered as Markov chains on a product space $\Gamma = S^W$ where $S$ is a finite set and $W$ is an infinite countable set of *sites*. The set $\Gamma$ is usually called *the configuration space*. Clearly, $\Gamma$ is a compact metrizable space with respect to the standard product discrete topology. Markov chains of this type were introduced originally by Stavskaya and Pyatetskii-Shapiro in [15] as a model for a neuron network and by Wasserstein in [16] as a model describing large systems of automata. Later, there have been extensive studies of various properties of PCA (see, for instance, [2], [10], [13], [12], as well as more recent works, like, for example, [1] or [12], and references there).

The large deviations for the temporal occupational measures (1.1) for this class of Markov chains have been investigated in [8] and [9], where some anomalous phenomena have been demonstrated. In particular, we have shown that for some natural examples of PCA, the lower bounds (1.4) given by the Donsker -Varadhan action functional $I(\nu)$ are not valid, while the upper bounds (1.3) are not optimal (although they are surely valid), and they can be replaced by a family of some alternative action functionals depending on the initial distribution of the Markov chain. On the other hand, we also demonstrated in Section 6 of [9] a class of PCA cases satisfying the bounds (1.3) and (1.4) (although under very strong additional conditions).

The author would like to emphasize that even though the upper bounds (1.3) are not optimal, they unquestionably hold for all typical PCA, as it has already been mentioned in the previous paragraph. For this reason, in order to delve deeper into the asymptotic behavior of the occupational measure $\xi_T$, it is crucial to comprehend the properties of the action functional $I(\nu)$. The purpose of the present paper is to investigate the Donsker-Varadhan action functional for infinite-dimensional Markov chains, with a focus on PCA. The main mathematical contribution of this investigation lies in demonstrating surprising features of the action functional under specific conditions for PCA, which provides some new insights into large deviations theory for infinite-dimensional systems.

To be more specific, it will be shown that for a wide family of PCA, the Donsker-Varadhan action functional $I(\nu)$ defined by the formula (1.2) obtains values different from infinity only on a rather restricted class of probability measures. The following reasoning will make the potential applications of this result evident. Suppose that some closed with respect to the weak topology subset K of $M(\Gamma)$ satisfies the condition $I(\nu) = \infty$ for all $\nu \in K$, then the inequality (1.3) yields the following fact:

$$\limsup_{T \to \infty} \frac{\log P_x \, (\xi_T \in K)}{T} = -\infty,$$

which means that the corresponding probability is asymptotically negligible even compared to the exponentially small probabilities. Therefore, it is important to prove the condition $I(\nu) = \infty$ for as wide a family of measures as is achievable in order to indicate negligible events. One possible example of such an approach is the formula (1.8) below. Moreover, the results presented here may also be applicable in proving the lower large deviation bounds (1.4) for a significantly wider class of probabilistic cellular automata than those considered in [9]. Indeed, when $I(\nu) = \infty$, the bounds in (1.4) follow trivially. Thus, the remaining challenge is to establish (1.4) for the restricted class of measures $\nu \in M(\Gamma)$ with $I(\nu) < \infty$, relying on the special properties derived in the present work. The author intends to pursue this direction in future research.

The paper is structured as follows: Section 2, which is broken into four subsections, contains the main results of this work. Namely, Subsection 2.1 provides the overall setup of the paper and formulates Lemma 2.1, which serves as the foundation for the key theorems that are formulated and proved in the following subsections (while the proof of this lemma is moved to Section 3, since it is rather technical and voluminous).

The main result, namely, **Theorem 1** is formulated and proved in Subsection 2.2 and indicates that

$$I(\nu) < \infty \tag{1.5}$$

only in the case when $\nu \in M(\Gamma)$ is a stationary measure of *the boundary process* associated to our Markov chain $X_t$. For the convenience of the reader, recall that Fölmer introduced the boundary process in the following way (see [10]). Let $\Phi$ be a finite subset of $W$, and let $\Phi^\star = W - \Phi$. Denote by $\mathcal{B}_{\Phi^\star}$ the sub-$\sigma$-algebra of $\mathcal{B}$ generated by the natural projection $\pi_{\Phi^\star} : \Gamma \to S^{\Phi^\star}$, which allows us to define *the spatial tail $\sigma$-algebra*

$$\widehat{\mathcal{B}} = \bigcap_{\Phi \subset W : \, \Phi \text{ is finite}} \mathcal{B}_{\Phi^\star} \tag{1.6}$$

Let $P(x, \cdot)$ be the transition probability kernel of on $(\Gamma, \mathcal{B})$. It turns out that under reasonably standard conditions on PCA (in particular, under the conditions (A1) - (A4) formulated in the next section), the function $P(\cdot, A)$ is $\widehat{\mathcal{B}}$-measurable for each $A \in \widehat{\mathcal{B}}$. Therefore, the transition probability kernel $P(x, \cdot)$ defines a Markov chain on $(\Gamma, \widehat{\mathcal{B}})$ which is called *the boundary process* associated to the Markov chain $X_t$. It will be shown in Subsection 2.2 that if $\nu \in M(\Gamma)$ satisfies the condition (1.5), then $\nu \in M(\Gamma)$ is a time-invariant measure of this boundary process, that is, for each $A \in \widehat{\mathcal{B}}$ we have

$$\int_\Gamma P(x, A)\nu(dx) = \nu(A). \tag{1.7}$$

Note that in this paper we will use the term *a time-invariant measure* to avoid confusion with space-shift invariant measures which are discussed in the next subsection.

Namely, in Subsection 2.3 we assume that $W = \mathbb{Z}^d$, that is, $\Gamma = S^{\mathbb{Z}^d}$, which allows us to consider *the space-shift transformations* $\theta_m: \Gamma \to \Gamma$ (defined in the usual way for any $m \in \mathbb{Z}^d$), and to introduce an additional assumption that the transition probability kernel of our Markov chain on $\Gamma$ is space-shift invariant. Under this assumption, **Theorem 2** states that if $\nu \in M(\Gamma)$ is a space-shift invariant measure, then either $I(\nu) = 0$ or $I(\nu) = \infty$. As a direct consequence of applying Theorem 2 to the bounds (1.3), the following is stated in Corollary 2.9. Denote by $M_{TI}(\Gamma)$ the set of all the time-invariant and by $M_{SSI}(\Gamma)$ the set of all the space-shift invariant probability Borel measures defined on $\Gamma$. Then for any open (with respect to the weak topology) neighborhood $U$ of $M_{TI}(\Gamma)$ it holds that

$$\underset{T \to \infty}{\limsup} \frac{\log P_x(\xi_T \in K)}{T} = -\infty, \tag{1.8}$$

where $K = \{\nu \in M_{SSI}(\Gamma): \nu \notin U\}$.

Finally, Subsection 2.4 draws attention to the fact that, under essentially the same assumptions, any measure satisfying constraint (1.5), while not required to be space-shift invariant, must satisfy some strong necessary conditions with respect to the space-shifts. To make this fact transparent, we will utilize a simplified set of sites, namely, $W = \{0,1,2,...\} = \mathbb{Z}^+$ in this subsection. Therefore, in this case, $\Gamma = S^{\mathbb{Z}^+}$, and one can use here the left shift transformations $\theta_n: \Gamma \to \Gamma$ defined for any integer $n > 0$. For any given measure $\mu \in M(\Gamma)$, these left shifts define the sequence of the left shifted measures $\mu_n$, by the formula $\mu_n(A) = \mu(\theta_n^{-1}(A))$. **Theorem 3** states that if there exists only one time-invariant measure $\nu_{TI} \in M(\Gamma)$ with respect to $P(x, \cdot)$, then for any measure $\mu$ satisfying condition (1.5), the sequence $\mu_n$ converges to $\nu_{TI}$ with respect to the weak topology.

To conclude, we remark that explicit examples of Markov chains satisfying the standing assumptions are omitted for brevity, as many such examples are available in the existing literature.

## 2. THE MAIN RESULTS

### 2.1. THE PCA SET-UP: PRELIMINARIES AND THE MAIN LEMMA

Let $S$ is a finite nonempty set and $W$ be an infinite countable set. Set $\Gamma = S^W$. As it was pointed out in Introduction, $\Gamma$ is a compact metrizable space with respect to the standard product discrete topology. We consider a transition probability kernel $P(x,\cdot)$ on $(\Gamma, \mathcal{B})$, where $\Gamma = S^W$ and $\mathcal{B}$ is the Borel $\sigma$-algebra of $\Gamma$. In order to provide a rigorous description of the relevant Markov chains, we will need some additional notations.

*Notations.*

a) For any subset $\Phi \subset W$ and any site $z \in \Phi$ define the natural projection $\pi_z: S^\Phi \to S$ by the formula $\pi_z(u) = u(z)$ for any $u \in S^\Phi$

b) Likewise, for any $\Phi \subset \Psi \subset W$ introduce the natural projection $\pi_\Phi: S^\Psi \to S^\Phi$, that is, $w = \pi_\Phi(u) \in S^\Phi$ is such that $\pi_z(w) = \pi_z(u)$ for any $z \in \Phi$.

c) For any $\Phi \subset W$ denote by $\mathcal{B}_\Phi$ the sub-$\sigma$-algebra of $\mathcal{B}$ generated by the projections $\pi_z : \Gamma \to S$, $z \in \Phi$. In particular, if $\Phi \subset W$ is finite, then the sub-$\sigma$-algebra $\mathcal{B}_\Phi$ is generated by the finite partition

$$\Lambda_\Phi = \{G_v : v \in S^\Phi \}, \qquad (2.1)$$

where for each $v \in S^\Phi$ the cylindrical set $G_v$ is defined in the following way:

$$G_v = \{ y \in \Gamma : \pi_\Phi(y) = v\}. \qquad (2.2)$$

**Remark 2.1** If $\Phi \subset W$ is not finite, still we can say that the sub-$\sigma$-algebra $\mathcal{B}_\Phi$ is generated by the family of the cylindrical sets $G_w = \{ y \in \Gamma : \pi_\Psi(y) = w\}$ defined for all the finite subsets $\Psi$ of $\Phi$ and for each $w \in S^\Psi$.

d) We will also need the following notations: for each $s \in S$ and for any $\Phi \subset W$, $z \in \Phi$, $v \in S^\Phi$, denote $H_{s,z} = \{ u \in \Gamma : \pi_z(u) = s\}$ and,

$$G_{v,z} = H_{\pi_z(v),z} = \{ u \in \Gamma : \pi_z(u) = \pi_z(v)\}. \qquad (2.3)$$

**Remark 2.2** Obviously, due to (2.2), for any $v \in S^\Phi$ one has $G_v = \bigcap_{z \in \Phi} G_{v,z}$.

Throughout this paper we adapt the following standard model of PCA (we use the terminology of [2]). First of all, we assume that the kernel $P(x, \cdot)$ is *synchronous*, that is, the following condition is satisfied

**(A1)** For all $x \in \Gamma$, any finite $\Phi \subset W$ and each $v \in S^\Phi$,

$$P(x, G_v) = \prod_{z \in \Phi} P(x, G_{v,z}),$$

where $G_v$ and $G_{v,z}$ has been defined in (2.2) and (2.3).

Moreover, we assume that the kernel is *local*, that is,

**(A2)** For any $z \in W$ there exist a finite set $N(z) \subset W$ and a *local transition probability kernel* $P_z: S^{N(z)} \times S \to [0,1]$ (that is to say, for each $h \in S^{N(z)}$ a probability distribution $P_z(h, \cdot)$ is defined on $S$), such that for all $x \in \Gamma$ and for each $s \in S$,

$$P(x, H_{s,z}) = P_z(\pi_{N(z)}(x), s).$$

where $H_{s,z}$ has been defined in the paragraph d) of the *Notations* above. Moreover, without loss of generality, we will assume that $z \in N(z)$.

**Remark 2.3** By A1 and A2, for all $x \in \Gamma$, for any finite $\Phi \subset W$ and for each $v \in S^\Phi$ one has,

$$P(x, G_v) = \prod_{z \in \Phi} P_z(\pi_{N(z)}(x), \pi_z(v)) \qquad (2.4)$$

where $G_v$ has been defined by the formula (2.2).

While the following two simple additional assumptions are not among the conditions customarily required by traditional PCA models, they are usually satisfied in the majority of standard examples.

**(A3)** For any $z \in W$ and for each $v \in S^{N(z)}$ and each $s \in S$ it holds that $P_z(v, s) > 0$

**(A4)** For any $z \in W$ the set $D(z) = \{z' \in W : z \in N(z')\}$ is finite.

**Remark 2.4** The assumption $z \in N(z)$ (see **A2**) yields $z \in D(z)$ for any $z \in W$.

Now we can formulate our basic result:

**Lemma 2.1** Let $P(x, \cdot)$ be a transition probability kernel on $(\Gamma, \mathcal{B})$, such that the conditions (A1) - (A4) are satisfied, and let $\Phi_n$, $n \geq 1$, be a sequence of finite subsets of $W$ such that $\Phi_n \subset \Phi_{n+1}$, $\bigcup_{n=1}^\infty \Phi_n = W$. Denote: $\Phi_n^* = W - \Phi_n$ for $n \geq 1$, and let $\mathcal{B}_{\Phi_n^*}$ be the sub-$\sigma$-algebra of $\mathcal{B}$ generated by the projections $\pi_z : \Gamma \to S$, $z \in \Phi_n^*$. For any given $v \in M(\Gamma)$ introduce the measure $v^P \in M(\Gamma)$ defined for any $A \in \mathcal{B}$ by the formula

$$v^P(A) = \int_\Gamma P(x, A) v(dx). \qquad (2.5)$$

If

$$I(\nu) < \infty, \tag{2.6}$$

then

$$\lim_{n \to \infty} \sup_{A \in \mathcal{B}_{\Phi_n^*}} |\nu^P(A) - \nu(A)| = 0 \tag{2.7}$$

*Proof.* We will prove Lemma 2.1 in Section 3.

## 2.2. THE BOUNDARY PROCESS

As in the introduction, let $\widehat{\mathcal{B}}$ be *spatial tail $\sigma$-algebra* of on $(\Gamma, \mathcal{B})$, that is,

$$\widehat{\mathcal{B}} = \bigcap_{\Phi \subset W: \Phi \text{ is finite}} \mathcal{B}_{\Phi^*} \tag{2.8}$$

Here, for any given finite subset $\Phi$ of $W$, we set: $\Phi^* = W - \Phi$, and $\mathcal{B}_{\Phi^*}$ is the sub-$\sigma$-algebra of $\mathcal{B}$ generated by the projection $\pi_z : \Gamma \to S$, $z \in \Phi^*$. If $P(x, \cdot)$ is a transition probability kernel on $(\Gamma, \mathcal{B})$ and the function $P(\cdot, A)$ is $\widehat{\mathcal{B}}$-measurable for each $A \in \widehat{\mathcal{B}}$. then the transition probability kernel $P(x, \cdot)$ defines a Markov chain on $(\Gamma, \widehat{\mathcal{B}})$ which is called *the boundary process* associated to the Markov chain $X_t$. We will now show that if $\nu \in M(\Gamma)$ satisfies the condition (1.5), then $\nu \in M(\Gamma)$ is a time-invariant measure of the boundary process.

More precisely, our main result for the general case is the following theorem.

**Theorem 1.** Let $P(x, \cdot)$ be a transition probability kernel on $(\Gamma, \mathcal{B})$, such that the conditions (A1) - (A4) are satisfied. If $A \in \widehat{\mathcal{B}}$, then the function $P(\cdot, A)$ is $\widehat{\mathcal{B}}$-measurable for any $A \in \widehat{\mathcal{B}}$. Moreover, for any measure $\nu \in M(\Gamma)$ such that

$$I(\nu) < \infty, \tag{2.9}$$

where the actional functional $I(\nu)$ has been introduced in (1.2), one has for any $A \in \widehat{\mathcal{B}}$

$$\int_\Gamma P(x, A)\nu(dx) = \nu(A) \tag{2.10}$$

*Proof.* The fact that $P(\cdot, A)$ is $\widehat{\mathcal{B}}$-measurable for any $A \in \widehat{\mathcal{B}}$ has been pointed out by Fölmer in [10] under similar conditions. For the convenience of the reader, it seems reasonable to provide a short proof of this fact. Let $\Phi$ be a finite subset of $W$. Denote

$$\Phi_0 = \bigcup_{y \in \Phi} D(y), \quad \Phi_0^* = W - \Phi_0. \tag{2.11}$$

By **A4** and (2.11) it is clear that $\Phi_0^* = \{z \in W : \Phi \cap N(z) = \emptyset\}$, and, therefore,

$$\bigcup_{z \in \Phi_0^*} N(z) \subset \Phi^*. \tag{2.12}$$

First, we will show, that for any $A \in \mathcal{B}_{\Phi_0^*}$ the function $P(\cdot, A)$ is $\mathcal{B}_{\Phi^*}$ measurable. By Remark 2.1, it is enough to proof this fact for all the cylindrical sets of the form $A =$ $\{u \in \Gamma: \pi_\Psi(u) = w\}$ defined for all the finite subsets $\Psi$ of $\Phi_0^*$ and for each $w \in S^\Psi$. By (2.4), for any such $A$,

$$P(x, A) = \prod_{z \in \Psi} P_z\big(\pi_{N(z)}(x), \pi_z(w)\big),  \qquad (2.13)$$

and, therefore, $P(\cdot, A)$ is $\mathcal{B}_{\Psi_1}$ measurable, where $\Psi_1 = \bigcup_{z \in \Psi} N(z)$. But $\Psi \subset \Phi_0^*$, and, thus, by (2.12), $\Psi_1 \subset \Phi^*$, which yields $\mathcal{B}_{\Psi_1} \subset \mathcal{B}_{\Phi^*}$. Therefore, $P(\cdot, A)$ is $\mathcal{B}_{\Phi^*}$ measurable, provided $A \in \mathcal{B}_{\Phi_0^*}$.

Next, let $A \in \widehat{\mathcal{B}}$, then, by (2.11) and (2.8), $A \in \mathcal{B}_{\Phi_0^*}$ for any finite subset $\Phi$ of $W$, since $\Phi_0$ is also finite. Thus, $P(\cdot, A)$ is $\mathcal{B}_{\Phi^*}$ measurable for any such $\Phi$, that is to say, $P(\cdot, A)$ is $\widehat{\mathcal{B}}$ measurable, as we wanted to show.

Now we will proof the formula (2.10). Let $\nu \in M(\Gamma)$ and $I(\nu) < \infty$. Since $W$ is an infinite countable set, it is clear that there exists a sequence of finite subsets $\Phi_n$, $n \geq 1$ of $W$ such that $\Phi_n \subset \Phi_{n+1}$, $\bigcup_{n=1}^{\infty} \Phi_n = W$. Thus, the conditions of Lemma 2.2 are satisfied. Clearly, by (2.8), $\widehat{\mathcal{B}} \subset \mathcal{B}_{\Phi_n^*}$ for any $n \geq 1$. Therefore, by (2.7)

$$\sup_{A \in \widehat{\mathcal{B}}} |\nu^P(A) - \nu(A)| \leq \lim_{n \to \infty} \sup_{A \in \mathcal{B}_{\Phi_n^*}} |\nu^P(A) - \nu(A)| = 0,$$

That is, $\nu^P(A) = \nu(A)$ for any $A \in \widehat{\mathcal{B}}$, which, actually, is the formula (2.10).

**Q.E.D**

## 2.3. THE CASE OF SHIFT- INVARIANT TRANSITION PROBABILITIES.

In this subsection we assume that $W = \mathbb{Z}^d$, that is, $\Gamma = S^{\mathbb{Z}^d}$, $d \geq 1$. Obviously, we can treat any $x \in \Gamma$ as the function $x: \mathbb{Z}^d \to S$ defined by the formula $x(z) = \pi_z(x)$, where $z \in \mathbb{Z}^d$. For any $m \in \mathbb{Z}^d$ introduce *the space-shift transformation* $\theta_m: \Gamma \to \Gamma$ acting in the usual way: $y = \theta_m(x)$ for any $x \in \Gamma$, where $y: \mathbb{Z}^d \to S$ is given by the formula $y(z) = \pi_{z+m}(x)$ for any $z \in \mathbb{Z}^d$. In the present subsection we assume that, in addition to the conditions (A1) -(A4), the transition probability is space-shift invariant, that is. the following condition holds:

**(A5)** For all $x \in \Gamma$, $m \in \mathbb{Z}^d$ and for any $A \in \mathcal{B}$,

$$P(x, A) = P\big(\theta_m(x), \theta_m(A)\big).$$

(Recall that $\mathcal{B}$ is the Borel $\sigma$-algebra of $\Gamma$)

**Remark 2.5**. It is convenient to introduce a more general definition of the space shift transformations $\theta_m: S^\Psi \to S^\Phi$ for any given $\Psi \subset \mathbb{Z}^d$ and for any $m \in \mathbb{Z}^d$, where $\Phi = \Psi - m$, in the following way: $v = \theta_m(w)$ for any $w \in S^\Psi$, where $v(z) = \pi_{z+m}(u)$ for any $z \in \Phi$. It is clear that for any the cylindrical set of the form

$$C_w = \{ u \in \Gamma : \pi_\Psi(u) = w \}$$

where $w \in S^\Psi$, we have

$$\theta_m(C_w) = \{ y \in \Gamma : \pi_\Phi(y) = \theta_m(w) \}, \tag{2.14}$$

where $\Phi = \Psi - m$.

**Remark 2.6.** One can reformulate the condition (A5) in terms of the local transition probability kernels $P_z : S^{N(z)} \times S \to [0,1]$. More precisely, it is clear that if the conditions (A1)-(A2) are satisfied, the condition (A5) is equivalent to the following assumptions: for all $z, m \in \mathbb{Z}^d$ we have

$$N(z) + m = N(z + m), \tag{2.15}$$

and the local transition probability kernels $P_z$ satisfy the condition:

$$P_{z+m}(\gamma, s) = P_z(\theta_m(\gamma), s)$$

for each $\gamma \in S^{N(z+m)}$, $s \in S$.

*Additional notations.*

We will say that the measure $v \in M(\Gamma)$ is *space-shift invariant* if for all $m \in \mathbb{Z}^d$ and for any $A \in \mathcal{B}$ it holds that

$$v(A) = v(\theta_m(A)). \tag{2.16}$$

Denote by $M_{SSI}(\Gamma)$ the set of all the space-shift invariant probability Borel measures defined on $\Gamma$. If $v \in M(\Gamma)$ is an invariant measure with respect to our Markov chain in the usual sense, that is, if for each $A \in \mathcal{B}$ one has

$$\int_\Gamma P(x, A) v(dx) = v(A), \tag{2.17}$$

we will say that $v \in M(\Gamma)$ is *a time-invariant measure* (to avoid confusion). Denote by $M_{TI}(\Gamma)$ the set of all the time-invariant probability Borel measures defined on $\Gamma$.

**Remark 2.7.** Clearly, if the set $M_{TI}(\Gamma)$ is a singleton and the condition (A5) is satisfied, then $M_{TI}(\Gamma) \subset M_{SSI}(\Gamma)$.

The main result of the present subsection is the following theorem.

**Theorem 2.** Let $P(x, \cdot)$ be a transition probability kernel on $(\Gamma, \mathcal{B})$, such that the conditions (A1) - (A5) are satisfied. If $v \in M_{SSI}(\Gamma)$ and $I(v) < \infty$, then $v \in M_{TI}(\Gamma)$.

*Proof.* For $n \geq 1$ denote: $\Phi_n = \{\, z = (z_1, z_2, \ldots, z_d) \in \mathbb{Z}^d : \sum_{i=1}^{n}|z_i| \leq n \,\}$. Then $\Phi_n$, $n \geq 1$, be a sequence of finite subsets of $W$ such that $\Phi_n \subset \Phi_{n+1}$, $\bigcup_{n=1}^{\infty} \Phi_n = \mathbb{Z}^d$. Thus, the conditions of Lemma 2.1 are satisfied, and, therefore, if $I(\nu) < \infty$, then by (2.7),

$$\lim_{n \to \infty} \sup_{A \in \mathcal{B}_{\Phi_n^*}} |\nu^P(A) - \nu(A)| = 0, \tag{2.18}$$

where $\Phi_n^* = \mathbb{Z}^d - \Phi_n$, the measure $\nu^P \in M(\Gamma)$ is defined in (2.5), and $\mathcal{B}_{\Phi_n^*}$ is the sub-$\sigma$-algebra of $\mathcal{B}$ generated by the projections $\pi_z : \Gamma \to S$, $z \in \Phi_n^*$.

Due to the condition (A5), if $\nu \in M_{SSI}(\Gamma)$, then $\nu^P \in M_{SSI}(\Gamma)$, and, therefore, for all $m \in \mathbb{Z}^d$ and for any $A \in \mathcal{B}$ we have

$$|\nu^P(A) - \nu(A)| = |\nu^P(\theta_m(A)) - \nu(\theta_m(A))| \tag{2.19}$$

Our goal is to show that for any $A \in \mathcal{B}$ one has

$$\nu^P(A) = \nu(A) \tag{2.20}$$

Since the $\sigma$-algebra $\mathcal{B}$ is generated by the family of all the cylindrical sets of the form

$$C_w = \{\, u \in \Gamma : \pi_\Psi(u) = w \,\} \tag{2.21}$$

defined for all the finite subsets $\Psi$ of $\mathbb{Z}^d$ and for each $w \in S^\Psi$, it is enough to prove the formula (2.21) for any $A = C_w$ defined by (2.21). For this purpose, observe that for any $n \geq 1$ and for any given finite subset $\Psi$ of $\mathbb{Z}^d$ one can find $m \in \mathbb{Z}^d$ such that

$$\Phi = \Psi - m \subset \Phi_n^*, \tag{2.22}$$

and, therefore, due to the formula (2.14), for each $w \in S^\Psi$ one has

$$\theta_m(C_w) = \{\, y \in \Gamma : \pi_\Phi(y) = \theta_m(w) \,\} \in \mathcal{B}_{\Phi_n^*}. \tag{2.23}$$

Therefore, by (2.19) and (2.23), we have

$$|\nu^P(C_w) - \nu(C_w)| = |\nu^P(\theta_m(C_w)) - \nu(\theta_m(A))| \tag{2.24}$$

$$\leq \sup_{A \in \mathcal{B}_{\Phi_n^*}} |\nu^P(A) - \nu(A)|,$$

which, together with (2.18), proves that (2.20) holds for any for any $A = C_w$ defined by (2.21).

**Q.E.D**

**Remark 2.8.** Observe that Theorem 2 states, actually, that if $\nu$ is a space-shift invariant measure, then either $I(\nu) = \infty$ or $I(\nu) = 0$.

The following result is a direct outcome of Theorem 2 combined with the Donsker and Varadhan result (1.3).

**Corollary 2.9.** Let $X_t$, $t \in \mathbb{Z}^+$, be a time homogeneous Markov chain defined on the phase space $\Gamma$ such that the corresponding transition probability kernel $P(x, \cdot)$ satisfies the conditions (A1) - (A5). Then for any open (with respect to the weak topology) neighborhood $U$ of $M_{TI}(\Gamma)$ it holds that

$$\limsup_{T \to \infty} \frac{\log P_x(\xi_T \in M_{SSI}(\Gamma) - U)}{T} = -\infty , \qquad (2.25)$$

where $\xi_T$ has been defined in (1.1), and $M_{SSI} - U = \{\nu \in M_{SSI}(\Gamma): \nu \notin U\}$

## 2.4. THE CONVERGENCE OF SHIFTED MEASURES TO THE UNIQUE TIME-INVARIANT MEASURE

Now we will consider the following simplified configuration to highlight some further necessary conditions that should be met by the measures $\mu \in M(\Gamma)$ such that $I(\mu) < \infty$.

In this subsection we assume that $W = \mathbb{Z}^+ = \{0,1,2,...\}$, that is, $\Gamma = S^{\mathbb{Z}^+}$. As in the previous subsection, we can treat any $x \in \Gamma$ as the function $x: \mathbb{Z}^+ \to S$ defined by the formula $x(z) = \pi_z(x)$, where $z \in \mathbb{Z}^+$. Now we can introduce the space shifts $\theta_n: \Gamma \to \Gamma$ in the following way: for a given integer $n > 0$ define $y = \theta_n(x)$ for each $x \in \Gamma$, where $y: \mathbb{Z}^+ \to S$ is given by the formula $y(z) = \pi_{z+n}(x)$ for any $z \in \mathbb{Z}^+$.

Note that, contrary to the previous subsection, the space-shift transformations $\theta_n: \Gamma \to \Gamma$ are not invertible. For any $A \in \mathcal{B}$ denote

$$\theta_n^{-1}(A) = \{y \in \Gamma: \theta_n(y) \in A\} \qquad (2.26)$$

Next, for any finite set of form $\Psi = \{a, a+1, ..., b\}$ introduce the set $\Psi(n) = \{a + n, a + 1 + n, ..., b + n\}$, where $a, b, n \in \mathbb{Z}^+$, and define the space shifts $\tilde{\theta}_n: S^{\Psi(n)} \to S^{\Psi}$ by the similar formula, that is, for a given integer $n > 0$ define $u = \tilde{\theta}_n(w)$ for each $w \in S^{\Psi(n)}$, where $u: \Psi \to S$ is given by the formula $u(z) = \pi_{z+n}(w)$ for any $z \in \Psi$.

**Remark 2.10.** It is obvious that the functions $\tilde{\theta}_n$ are invertible and that for any integer $n > 0$ and for each $x \in \Gamma$ one has

$$\pi_\Psi(\theta_n(x)) = \tilde{\theta}_n\left(\pi_{\Psi(n)}(x)\right) \qquad (2.27)$$

Moreover, for a given $u \in S^\Psi$ let $G_u = \{y \in \Gamma : \pi_\Psi(y) = u\}$ be the corresponding cylindrical set defined as in (2.2). It follows immediately from (2.26) and (2.27) that

$$\theta_n^{-1}(G_u) = G_w, \qquad (2.28)$$

where $G_w = \{y \in \Gamma : \pi_{\Psi(n)}(y) = w\}$ and $w \in S^{\Psi(n)}$ is such that $u = \tilde{\theta}_n(w)$.

In the present subsection we assume, in addition to the conditions (A1)-(A4), that the local transition probability kernels are space-shift invariant, that is, the following condition holds:

**(A6)** For all $z \in \mathbb{Z}^+$ the set $N(z)$ is of the form

$$N(z) = \{z, z+1, \ldots, z+r_0\} \qquad (2.29)$$

for some given positive integer $r_0$ independent of $z$, and the local transition probability kernels $P_z$ satisfy the condition:

$$P_{z+n}(w, s) = P_z(\tilde{\theta}_n(w), s)$$

for each $w \in S^{N(z+n)}$, each $s \in S$ and any integer $n > 0$.

**Remark 2.11.** Clearly, if $\Psi = N(z)$, then $\Psi(n) = N(z+n)$, and, therefore, $\tilde{\theta}_n(w)$ is well defined. Moreover, the equality (2.27) takes the following form

$$\pi_{N(z)}(\theta_n(x)) = \tilde{\theta}_n\left(\pi_{N(z+n)}(x)\right)$$

for any integer $n > 0$, $z \geq 0$, and for each $x \in \Gamma$.

**Notation.** For any given measure $\mu \in M(\Gamma)$, the space shifts naturally define the sequence of measures $\mu_n$, $n \geq 1$, by the formula

$$\mu_n(A) = \mu(\theta_n^{-1}(A)) \qquad (2.30)$$

for any $A \in \mathcal{B}$.

**Remark 2.12.** Observe that the definition (2.30) is equivalent to the statement that the following formula

$$\int_\Gamma f(x) \mu_n(dx) = \int_\Gamma f(\theta_n(x)) \mu(dx)$$

holds for any bounded $\mathcal{B}$-measurable function $f$.

It turns out that, under the assumptions of this subsection, the following theorem holds.

**Theorem 3.** Assume that the transition probability kernel $P(x, \cdot)$ on $(\Gamma, \mathcal{B})$ is such that the conditions (A1) - (A4) and (A6) are satisfied, and, additionally, assume that there exists only one time-invariant measure $\nu_{TI} \in M(\Gamma)$ with respect to $P(x, \cdot)$. Let $\mu \in M(\Gamma)$ satisfies the condition: $I(\mu) < \infty$. Then the sequence of measures $\mu_n$, defined by (2.30), converges weakly to $\nu_{TI}$.

**Remark 2.13.** Recall that constructive sufficient conditions for the uniqueness of the time-invariant measure for PCA were already established in the classical work [16] by Wasserstein. Under assumption **(A6)**, these conditions take on an especially simple form.

*Proof of Theorem 3.* Again, our proof will be based on Lemma 2.1. Indeed, for any integer $n > 0$ define $\Phi_n^* = \{n, n+1, n+2, ...\}$, that is,

$$\Phi_n^* = \mathbb{Z}^+ - \Phi_n,$$

where $\Phi_n = \{0, 1, 2, ..., n-1\}$. Clearly, the assumptions of Lemma 2.1 are satisfied, and, therefore,

$$\lim_{n \to \infty} \sup_{A \in \mathcal{B}_{\Phi_n^*}} |\mu^P(A) - \mu(A)| = 0, \tag{2.33}$$

where $\mathcal{B}_{\Phi_n^*}$ be the sub-$\sigma$-algebra of $\mathcal{B}$ generated by the projections $\pi_z : \Gamma \to S, z \in \Phi_n^*$.

Next, we will show, using (2.33), that for any cylindrical sets of the form $G_u = \{y \in \Gamma : \pi_\Psi(y) = u\}$ defined in (2.2), such that $u \in S^\Psi$ and $\Psi = \{0, 1, ..., b\}$, $b \in \mathbb{Z}^+$, one has

$$\lim_{n \to \infty} |\mu_n^P(G_u) - \mu_n(G_u)| = 0. \tag{2.34}$$

We will prove the assertion (2.34) through the formulae (2.35) to (2.39).

Indeed, for a given $u \in S^\Psi$ and an integer $n > 0$ let $w \in S^{\Psi(n)}$ be such that $u = \tilde{\theta}_n(w)$, where $\tilde{\theta}_n$ and $\Psi(n) = \{n, n+1, ..., n+b\}$ are defined as in the beginning of this subsection. Then, by (2.30) and (2.28),

$$\mu_n(G_u) = \mu(\theta_n^{-1}(G_u)) = \mu(G_w) \tag{2.35}$$

where $G_w = \{y \in \Gamma : \pi_{\Psi(n)}(y) = w\}$. Since $\Psi(n) \subset \Phi_n^*$, it is clear that

$$G_w \in \mathcal{B}_{\Phi_n^*}. \tag{2.36}$$

On the other hand, , substituting $w$ instead of $v$ in (2.4) and applying (A6), then taking into account Remark 2.11 and the definition of $w$, or, more precisely, the fact that $\pi_z(u) = \pi_{z+n}(w)$ for any $z \in \Psi$, and applying again (2.4), we obtain for any $x \in \Gamma$ the following convenient formula (similar to the condition (A5) of the previous subsection),

$$P(x, G_w) = \prod_{z \in \Psi(n)} P_z(\pi_{N(z)}(x), \pi_z(w)) = \tag{2.37}$$

$$\prod_{n \le z \le n+b} P_z(\pi_{N(z)}(x), \pi_z(w)) = \prod_{0 \le z \le b} P_{z+n}(\pi_{N(z+n)}(x), \pi_{z+n}(w)) =$$

$$\prod_{0 \le z \le b} P_z\left(\tilde{\theta}_n\left(\pi_{N(z+n)}(x)\right), \pi_{z+n}(w)\right) = \prod_{z \in \Psi} P_z(\pi_{N(z)}(\theta_n(x)), \pi_z(u)) = P(\theta_n(x), G_u)$$

Now, by the definition (2.5), together with Remark 2.12, and by the formula (2.37),

$$\mu_n^P(G_u) = \int_\Gamma P(x, G_u)\mu_n(dx) = \int_\Gamma P(\theta_n(x), G_u)\mu(dx) = \quad (2.38)$$

$$\int_\Gamma P(x, G_w)\mu(dx) = \mu^P(G_w)$$

Bringing together (2.38), (2.35) and (2.36), we obtain

$$|\mu_n^P(G_u) - \mu_n(G_u)| = |\mu^P(G_w) - \mu(G_w)| \leq \sup_{A \in \mathcal{B}_{\Phi_n^*}} |\mu^P(A) - \mu(A)|, \quad (2.39)$$

which, together with (2.33), proves (2.34).

Next, since $M(\Gamma)$ is compact with respect to the weak convergence topology, it is enough to prove that if a subsequence of $\mu_n$ converges weakly to some measure $\nu \in M(\Gamma)$, then

$$\nu = \nu_{TI} \quad (2.40)$$

More specifically, let $\mu_{n_m}$ be a subsequence of $\mu_n$ such that $\mu_{n_m}$ converges weakly to a measure $\nu \in M(\Gamma)$. Our purpose is to prove that for any $A \in \mathcal{B}$,

$$\nu^P(A) = \nu(A) \quad (2.41)$$

where $\nu^P(A)$ is defined in (2.5), which, together with the assumption concerning the uniqueness of the time-invariant measure, yields (2.31).

Observe that since the $\sigma$-algebra $\mathcal{B}$ is generated by the family of all the cylindrical sets of the form $G_u = \{y \in \Gamma : \pi_\Psi(y) = u\}$ defined in (2.2), such that $u \in S^\Psi$, where $\Psi = \{0, 1, \ldots, b\}$, $b \in \mathbb{Z}^+$, it is enough to prove (2.41) for any $A = G_u$. Therefore, in order to complete the proof, it remains to show that for each $u \in S^\Psi$ we have

$$\nu^P(G_u) = \nu(G_u) \quad (2.42)$$

Indeed, by (2.34),

$$\lim_{m \to \infty} |\mu_{n_m}^P(G_u) - \mu_{n_m}(G_u)| = 0 \ . \quad (2.43)$$

On the other hand, the indicator function $I_{G_u}(x)$ is continuous with respect to our topology, and, therefore,

$$\lim_{m \to \infty} \mu_{n_m}(G_u) = \lim_{m \to \infty} \int_\Gamma I_{G_u}(x)\mu_{n_m}(dx) = \int_\Gamma I_{G_u}(x)\nu(dx) = \nu(dx) \quad (2.44)$$

Similarly, by (2.4), for any $x \in \Gamma$, the function $P(x, G_u)$ satisfies the following equality

$$P(x, G_u) = \prod_{z \in \Psi} P_z\bigl(\pi_{N(z)}(x), \pi_z(u)\bigr)$$

and, since $\Psi$ is a finite subset of $\mathbb{Z}^+$, it is continuous in our topology with respect to $x \in \Gamma$. Therefore,

$$\lim_{m \to \infty} \mu_{n_m}^P(G_u) = \lim_{m \to \infty} \int_\Gamma P(x, G_u)\mu_{n_m}(dx) = \int_\Gamma P(x, G_u)\nu(dx) = \nu^P(G_u). \quad (2.45)$$

Clearly, the formulae (2.43) - (2.45) yield (2.42), hence completing the proof.

**Q.E.D**

**Corollary 2.14.** Under the conditions of Theorem 3, if the measure $\mu \in M(\Gamma)$ is space-shift invariant and $I(\mu) < \infty$, then $\mu$ it is time-invariant.

*Proof.* The statement follows immediately from the definition of space-invariant measure, that is, from the fact that $\mu(A) = \mu(\theta_n^{-1}(A))$ for any $A \in \mathcal{B}$, and for any integer $n \geq 1$, together with the definition (2.30).

## 3. PROOF OF LEMMA 2.1

The proof of Lemma 2.1 is separated into three stages. In Subsection 3.1, we prove Proposition 3.4 for a more abstract set-up. Following that, in Subsection 3.2, we prove Proposition 3.5, asserting that the setup of Subsection 3.1 is correct under our primary assumptions (A1-A4). As a result, given these assumptions, Proposition 3.4 holds and takes the form of Corollary 3.6, which, when paired with some standard properties of the action functional, yields the main result of Subsection 3.2, namely, Corollary 3.7. Finally, in Subsection 3.3, we conclude the proof of Lemma 2.1.

### 3.1. THE AUXILARY ESTIMATE

The goal of this subsection is to prove Proposition 3.4 below which serves as the key part of the proof of Lemma 2.1. We shall establish Proposition 3.4 for a much more general setup, since the author conjectured that it may be possible to prove some theorems like Lemma 2.1 not only for PCA but also for more general classes of infinite-dimensional systems. Because of this, this subsection will employ the more general and self-contained notation system that follows.

*Assumptions and notations of Subsection 3.1.*

a) Let $\Gamma = H \times \Gamma_0$ where $H$ is a finite set and $(\Gamma_0, \mathcal{B}_0)$ is a measurable space. Denote by $\pi_H : \Gamma \to H$ and by $\pi_0 : \Gamma \to \Gamma_0$ the natural projections from $\Gamma$ to $H$ and to $\Gamma_0$, respectively. Introduce the following $\sigma$-algebra of subsets of $\Gamma$,

$$\mathcal{B}'_0 = \{\pi_0^{-1}(B): B \in \mathcal{B}_0\} \quad (3.1)$$

and define the following finite partition of $\Gamma$,

$$\Lambda = \{A_h = \pi_H^{-1}(h): h \in H\} \quad (3.2)$$

Denote by $\mathcal{B}$ the $\sigma$-algebra generated by $\mathcal{B}'_0$ and $\Lambda$.

b) Next, define on the set $\Gamma \times \Gamma$ the left and the right projections $\tilde{\pi}_L$ and $\tilde{\pi}_R$, that is to say, $\tilde{\pi}_L(u,v) = u$ and $\tilde{\pi}_R(u,v) = v$ for any $(u,v) \in \Gamma \times \Gamma$. Denote by $\mathcal{A}$ the $\sigma$-algebra of subsets of $\Gamma \times \Gamma$ generated by the projections $\tilde{\pi}_L$ and $\tilde{\pi}_R$ with respect to $\mathcal{B}$, that is, $\mathcal{A} = \mathcal{B} \otimes \mathcal{B}$, and denote by $M(\Gamma \times \Gamma)$ the set of all the probability measures defined on $(\Gamma \times \Gamma, \mathcal{A})$.

c) Let $\mu \in M(\Gamma \times \Gamma)$ and let $\mu_L$ be the left marginal of $\mu$. For a given transition probability kernel $P(x, \cdot)$ on $(\Gamma, \mathcal{B})$, denote by $\mu^P \in M(\Gamma \times \Gamma)$ the measure defined by the formula

$$\mu^P(A \times B) = \int_A P(x, B) \mu_L(dx) \qquad (3.3)$$

for any $A, B \in \mathcal{B}$. Let $D(\mu \parallel \mu^P)$ be the relative entropy of $\mu$ with respect to $\mu^P$. Recall that if $\mu \ll \mu^P$, then

$$D(\mu \parallel \mu^P) = \int_{\Gamma \times \Gamma} \ln\left(\frac{d\mu}{d\mu^P}\right) d\mu. \qquad (3.4)$$

If $\mu$ is not absolutely continuous with respect to $\mu^P$, define, as usual,

$$D(\mu \parallel \mu^P) = \infty \qquad (3.5)$$

Throughout this subsection we will assume that a given transition probability kernel $P(x, \cdot)$ on $(\Gamma, \mathcal{B})$ satisfies the following condition.

**Assumption 3.1.** There exist a subpartition $\Lambda_1$ of $\Lambda$ and a sub-$\sigma$-algebra $\mathcal{B}'_1$ of $\mathcal{B}'_0$ and transition probability kernels $P_H: H \times \Lambda_1 \to [0,1]$ and $P_0: \Gamma_0 \times \mathcal{B}'_1 \to [0,1]$ such that for any $x \in \Gamma$ and for each $B \in \Lambda_1, C \in \mathcal{B}'_1$,

$$P(x, B \cap C) = P_H(\pi_H(x), B) P_0(\pi_0(x), C) \qquad (3.6)$$

Furthermore, we assume that for each $h \in H$, $B \in \Lambda_1$,

$$P_H(h, B) > 0 \qquad (3.7)$$

*Additional notations*

d) In the following Propositions 3.2 and 3.4 below we will need some additional notations. Introduce the sub-$\sigma$-algebra $\mathcal{A}_0 = \mathcal{B}'_0 \otimes \mathcal{B}'_1$ of $\mathcal{A}$ generated by the semi-algebra of the sets of form $D \times C \in \mathcal{A}$ such that $D \in \mathcal{B}'_0, C \in \mathcal{B}'_1$, and define the finite partition $\Delta = \Lambda \times \Lambda_1$ of $\Gamma \times \Gamma$, where $\Lambda_1$ and $\mathcal{B}'_1$ have been introduced in Assumption 3.1.

e) Finally, for a given $\mu \in M(\Gamma \times \Gamma)$, denote by $D_\Delta(\mu \parallel \mu^P)$ the relative entropy of $\mu$ with respect to $\mu^P$ corresponding to the partition $\Delta$, that is,

$$D_\Delta(\mu \parallel \mu^P) = \sum_{A \times B \in \Delta} \mu(A \times B) \ln\left(\frac{\mu(A \times B)}{\mu^P(A \times B)}\right), \qquad (3.8)$$

where we adapt the usual conventions concerning the relative entropy: $0\ln 0 = 0$

and $0\ln\frac{0}{0} = 0$ ( The reader will convince himself that $D_\Delta(\mu \parallel \mu^P)$ is well defined, due to Remark 3.3 below).

Observe that Assumption 3.1 yields the following simple results that will be useful in the further calculations:

**Proposition 3.2**

**a)** For any $x \in \Gamma$ and for any $C \in \mathcal{B}'_1$, one has

$$P(x, C) = P_0(\pi_0(x), C) \tag{3.9}$$

**b)** For each $h \in H$, $B \in \Lambda_1$, and for any $x \in A_h \in \Lambda$,

$$P(x, B) = P_H(h, B) \tag{3.10}$$

**c)** For any $\mu \in M(\Gamma \times \Gamma)$, for each set $A_h \times B \in \Delta$ and

for all $G \in \mathcal{A}_0$ we have,

$$\mu^P\big((A_h \times B) \cap G\big) = P_H(h, B)\mu^P\big((A_h \times \Gamma) \cap G\big) \tag{3.11}$$

In particular,

$$\mu^P(A_h \times B) = P_H(h, B)\mu_L(A_h), \tag{3.12}$$

(recall that the sets $A_h$ have been introduced in (3.2)).

*Proof.* a) Since $P_H(\pi_H(x), \Gamma) = 1$, it follows immediately from (3.6).

b) Similarly, since $\pi_H(x) = h$ and $P_0(\pi_0(x), \Gamma) = 1$.

c) Recall that $\mathcal{A}_0 = \mathcal{B}'_0 \otimes \mathcal{B}'_1$. Therefore, it is sufficient to prove (3.11) for all the sets $G \in \mathcal{A}_0$ of the form $G = D \times C$, where $D \in \mathcal{B}'_0$, $C \in \mathcal{B}'_1$. However, for such sets one can derive the formula (3.11) by the following direct calculation:

$$\begin{aligned}
\mu^P\big((A_h \times B) \cap G\big) &= \mu^P\big((A_h \times B) \cap (D \times C)\big) \\
&= \mu^P\big((A_h \cap D) \times (B \cap C)\big) = \int_{A_h \cap D} P(x, B \cap C)\mu_L(dx) \\
&= \int_{A_h \cap D} P_H(h, B)P_0(\pi_0(x), C)\mu_L(dx) =
\end{aligned} \tag{3.13}$$

$$\begin{aligned}
&= P_H(h, B) \int_{A_h \cap D} P_0(\pi_0(x), C)\mu_L(dx) \\
&= P_H(h, B) \int_{A_h \cap D} P(x, C)\mu_L(dx) = P_H(h, B)\mu^P\big((A_h \cap D) \times C\big) = \\
&= P_H(h, B)\mu^P\big((A_h \times \Gamma) \cap (D \times C)\big) = \\
&= P_H(h, B)\mu^P\big((A_h \times \Gamma) \cap G\big),
\end{aligned}$$

using (3.3), (3.6), (3.9) and again (3.3).

Finally, if $G = \Gamma \times \Gamma$, the formula (3.11) yields (3.12). Indeed, using (3.3) and the fact that $P(x, \Gamma) = 1$, one has

$$\mu^P(A_h \times B) = P_H(h, B)\mu^P(A_h \times \Gamma) = P_H(h, B)\mu_L(A_h).$$

**Q.E.D**

**Remark 3.3** By (3.12) and (3.7), if $\mu^P(A_h \times B) = 0$, then $\mu(A_h \times B) = 0$.
Therefore, $D_\Delta(\mu \| \mu^P)$ is well defined, and moreover, by (3.12), the definition (3.8) takes the following form

$$D_\Delta(\mu \| \mu^P) = \sum_{A_h \times B \in \Delta} \mu(A_h \times B) \ln\left(\frac{\mu(A \times B)}{\mu_L(A_h) P_H(h,B)}\right). \tag{3.14}$$

The next proposition is the main result of this subsection.

**Proposition 3.4** Let $\mu \in M(\Gamma \times \Gamma)$ be such that

$$D(\mu \| \mu^P) < \infty. \tag{3.15}$$

Then, under Assumption 3.1,

$$\sup_{G \in \mathcal{A}_0} |\mu^P(G) - \mu(G)| \leq \frac{1}{\sqrt{2}} \left( D(\mu \| \mu^P) - D_\Delta(\mu \| \mu^P) \right)^{\frac{1}{2}}. \tag{3.16}$$

*Proof of Proposition 3.4.* Let $\mu \in M(\Gamma \times \Gamma)$ be such that condition (3.15) is satisfied. For each $A_h \times B \in \Delta$ such that $\mu(A_h \times B) > 0$ define the conditional probability measure $\mu_{h,B}$ with respect to the event $A_h \times B$ on the measurable space $(\Gamma \times \Gamma, \mathcal{A}_0)$ in the standard way, that is, for any $G \in \mathcal{A}_0$ define

$$\mu_{h,B}(G) = \frac{1}{\mu(A_h \times B)} \mu\big((A_h \times B) \cap G\big). \tag{3.17}$$

Similarly, for each $A_h \times B \in \Delta$ such that $\mu(A_h \times B) > 0$ ( and, therefore, such that $\mu_L(A_h) > 0$) define the conditional probability measure $\mu^P_{h,B}$ on $(\Gamma \times \Gamma, \mathcal{A}_0)$ by the following formula, using (3.12),

$$\begin{aligned}\mu^P_{h,B}(G) &= \frac{1}{\mu^P(A_h \times B)} \mu^P\big((A_h \times B) \cap G\big) = \\ &= \frac{1}{\mu_L(A_h) P_H(h,B)} \mu^P\big((A_h \times B) \cap G\big),\end{aligned} \tag{3.18}$$

for any $G \in \mathcal{A}_0$. Recall that, due to (3.15), $\mu$ is absolutely continuous with respect to $\mu^P$, and, hence, if $\mu^P_{h,B}(G) = 0$, then $\mu_{h,B}(G) = 0$. Therefore, $\mu_{h,B}$ is absolutely continuous with respect to $\mu^P_{h,B}$. Denote

$$\rho_{h,B} = \frac{d\mu_{h,B}}{d\mu^P_{h,B}} \tag{3.19}$$

(For the sake of neatness, note that if $\mu(A_h \times B) = 0$, $\mu^P(A_h \times B) \neq 0$, then $\mu^P_{h,B}$ is well defined and we can define $\mu_{h,B} = \mu^P_{h,B}$. If $(A_h \times B) = \mu^P(A_h \times B) = 0$, define $\mu_{h,B} =$ $\mu_{h,B}^P = \mu$. In both situations, we can write $\rho_{h,B} = 1$.) Now we can define is the relative entropy of $\mu_{h,B}$ with respect to $\mu_{h,B}^P$

$$D(\mu_{h,B} \parallel \mu_{h,B}^P) = \int_{\Gamma \times \Gamma} \ln(\rho_{h,B}) \, d\mu_{h,B}. \tag{3.20}$$

(If $\mu(A_h \times B) = 0$, then, according to our conventions, $D(\mu_{h,B} \parallel \mu_{h,B}^P) = 0$.)

Observe that, due to (3.17) and (3.20), for each $A_h \times B \in \Delta$ the following formula holds true

$$\int_{\Gamma \times \Gamma} \chi_{A_h \times B} \ln(\rho_{h,B}) \, d\mu = \mu(A_h \times B) \, D(\mu_{h,B} \parallel \mu_{h,B}^P), \tag{3.21}$$

(here $\chi_{A_h \times B}$ is the corresponding indicator function).

Next, let $\mathcal{A}_1$ be the sub-$\sigma$-algebra of $\mathcal{A}$ generated by the partition $\Delta$ and the $\sigma$-algebra $\mathcal{A}_0$. Clearly, $\mathcal{A}_1$ consists of all the finite unions of the sets of the form $(A_h \times B) \cap G$, where $G \in \mathcal{A}_0, h \in H, B \in \Lambda_1$. Therefore, due to (3.17), (3.18) and (3.19), for each $A_h \times B \in \Delta$ such that $\mu(A_h \times B) > 0$, the Radon-Nikodym derivative of $\mu$ with respect to $\mu^P$ restricted to $\mathcal{A}_1$ could be given by the formula

$$\rho_{\mathcal{A}_1}(x, y) = \frac{\mu(A_h \times B)}{\mu_L(A_h) P_H(h, B)} \rho_{h,B}(x, y), \tag{3.22}$$

provided $(x, y) \in A_h \times B$.

Consider the relative entropy $D_{\mathcal{A}_1}(\mu \parallel \mu^P)$ of $\mu$ with respect to $\mu^P$ on the $\sigma$-algebra $\mathcal{A}_1$, that is

$$D_{\mathcal{A}_1}(\mu \parallel \mu^P) = \int_{\Gamma \times \Gamma} \ln(\rho_{\mathcal{A}_1}) \, d\mu. \tag{3.23}$$

It is a well-known fact that,

$$D_{\mathcal{A}_1}(\mu \parallel \mu^P) \leq D(\mu \parallel \mu^P) \tag{3.24}$$

(see, for instance, Corollary 5.2.2 of [11]). On the other hand, by (3.23), (3.22), (3.21) and (3.13), one has,

$$\begin{aligned} D_{\mathcal{A}_1}(\mu \parallel \mu^P) &= \sum_{A_h \times B \in \Delta} \int_{A_h \times B} \ln(\rho_{\mathcal{A}_1}) \, d\mu \\ &= \sum_{A_h \times B \in \Delta: \mu(A_h \times B) \neq 0} \int_{A_h \times B} \ln\left(\frac{\mu(A_h \times B)}{\mu_L(A_h) P_H(h, B)} \rho_{h,B}\right) d\mu \\ &= \sum_{A_h \times B \in \Delta: \mu(A_h \times B) \neq 0} \mu(A_h \times B) \ln\left(\frac{\mu(A_h \times B)}{\mu_L(A_h) P_H(h, B)}\right) \\ &\quad + \sum_{A_h \times B \in \Delta} \int_{\Gamma \times \Gamma} \chi_{A_h \times B} \ln(\rho_{h,B}) \, d\mu \\ &= D_\Delta(\mu \parallel \mu^P) + \sum_{A_h \times B \in \Delta} \mu(A_h \times B) D(\mu_{h,B} \parallel \mu_{h,B}^P). \end{aligned} \tag{3.25}$$

Observe that the formula (3.25) could be considered as some modification of the well-known chain rule for the relative entropy (see, for instance, [3], Theorem C.3.1), however, in our case it is easier to derive it directly from the definitions.

By (3.25) and (3.24),

$$\sum_{A_h \times B \in \Delta} \mu(A_h \times B) D(\mu_{h,B} \| \mu^P_{h,B}) \leq D(\mu \| \mu^P) - D_\Delta(\mu \| \mu^P) \quad (3.26)$$

Next, for each $h \in H$, $B \in \Lambda_1$ we can use the following fundamental inequality, involving the relative entropy and the variational distance between two given probability measures,

$$\sup_{G \in \mathcal{A}_0} |\mu_{h,B}(G) - \mu^P_{h,B}(G)| \leq \sqrt{\frac{D(\mu_{h,B} \| \mu^P_{h,B})}{2}} \quad (3.27)$$

(see [11], Lemmas 5.2.7 and 5.2.8).

Let $G \in \mathcal{A}_0$, then, combining the Cauchy–Schwarz inequality and the estimates (3.27) and (3.26), we obtain the following estimate

$$\begin{aligned}
&\sum_{A_h \times B \in \Delta} \mu(A_h \times B) \, |\mu_{h,B}(G) - \mu^P_{h,B}(G)| \\
&\leq \frac{1}{\sqrt{2}} \sum_{A_h \times B \in \Delta} (\mu(A_h \times B))^{\frac{1}{2}} \left(\mu(A_h \times B) D(\mu_{h,B} \| \mu^P_{h,B})\right)^{\frac{1}{2}} \\
&\leq \frac{1}{\sqrt{2}} \left(\sum_{A_h \times B \in \Delta} \mu(A_h \times B)\right)^{\frac{1}{2}} \left(\sum_{A_h \times B \in \Delta} \mu(A_h \times B) D(\mu_{h,B} \| \mu^P_{h,B})\right)^{\frac{1}{2}} \\
&= \frac{1}{\sqrt{2}} \left(\sum_{A_h \times B \in \Delta} \mu(A_h \times B) D(\mu_{h,B} \| \mu^P_{h,B})\right)^{\frac{1}{2}} \leq \frac{1}{\sqrt{2}} \left(D(\mu \| \mu^P) - D_\Delta(\mu \| \mu^P)\right)^{\frac{1}{2}}
\end{aligned} \quad (3.28)$$

Next, let us consider in more details the measure $\mu^P$ restricted to the sub-$\sigma$-algebra $\mathcal{A}_0$ and the measures $\mu^P_{h,B}$. By (3.18) and (3.11), for each $h \in H$ such that $\mu_L(A_h) > 0$, and for each $B \in \Lambda_1$ and any $G \in \mathcal{A}_0$,

$$\begin{aligned}
\mu^P_{h,B}(G) &= \frac{1}{\mu_L(A_h) P_H(h,B)} \mu^P\big((A_h \times B) \cap G\big) \\
&= \frac{1}{\mu_L(A_h)} \mu^P\big((A_h \times \Gamma) \cap G\big),
\end{aligned} \quad (3.29)$$

and, therefore, $\mu^P_{h,B}(G)$ is, actually, independent of $B \in \Lambda_1$. Thus, by (3.29), for any $G \in \mathcal{A}_0$, $B \in \Lambda_1$,

$$\mu^P(G) = \sum_{h \in H} \mu^P\big((A_h \times \Gamma) \cap G\big) = \sum_{h \in H} \mu_L(A_h) \mu^P_{h,B}(G). \quad (3.30)$$

Moreover, using the fact that $\mu^P_{h,B}(G)$ is independent of $B \in \Lambda_1$ and that

$$\mu_L(A_h) = \sum_{B \in \Lambda_1} \mu(A_h \times B),$$

one can rewrite (3.30) in the form

$$\mu^P(G) = \sum_{h \in H} \sum_{B \in \Lambda_1} \mu(A_h \times B) \mu_{h,B}^P(G). \tag{3.31}$$

On the other hand, clearly, (3.17) implies that

$$\mu(G) = \sum_{h \in H} \sum_{B \in \Lambda_1} \mu(A_h \times B) \mu_{h,B}(G). \tag{3.32}$$

Finally, by (3.31), (3.32) and (3.28),

$$|\mu^P(G) - \mu(G)| = \left| \sum_{h \in H} \sum_{B \in \Lambda_1} \mu(A_h \times B) \mu_{h,B}^P(G) - \sum_{h \in H} \sum_{B \in \Lambda_1} \mu(A_h \times B) \mu_{h,B}(G) \right|$$

$$\leq \sum_{A_h \times B \in \Delta} \mu(A_h \times B) |\mu_{h,B}(G) - \mu_{h,B}^P(G)| \leq \frac{1}{\sqrt{2}} \left( D(\mu \| \mu^P) - D_\Delta(\mu \| \mu^P) \right)^{\frac{1}{2}},$$

proving the estimate (3.16).
**Q.E.D**

### 3.2. RETURNING TO THE PCA SET-UP

In this subsection we will return to our original set-up introduced in Section 2. More precisely, we will show how our original PCA set-up fits into the more general framework developed in the previous subsection.

Assume that the conditions **(A1)** -**(A4)** of Section 2 are satisfied. Recall that $\Gamma = S^W$, where $S$ is a finite nonempty set and $W$ is an infinite countable set, and that $\mathcal{B}$ is the Borel $\sigma$-algebra of $\Gamma$. Throughout all this subsection, $\Phi$ is a given finite nonempty subset of $W$. The purpose of the present subsection is to prove Corollary 3.7 below. But, first, introduce the following notations.

*Notations of Subsection 3.2.* Denote

$$N(\Phi) = \bigcup_{z \in \Phi} N(z) \tag{3.33}$$

and

$$D(\Phi) = \bigcup_{y \in N(\Phi)} D(y) \tag{3.34}$$

Observe that the sets $D(\Phi)$ and $N(\Phi)$ are finite, due to the assumptions **(A2)** and **(A4)**, and that $\Phi \subset N(\Phi) \subset D(\Phi)$, due to Remark 2.4.

**Remark 3.5.** One can easily see that if $z \notin D(\Phi)$, then $N(z) \cap N(\Phi) = \emptyset$.

Indeed, if $y \in N(z)$, then, by the definition, $z \in D(y)$. On the other hand, if $y \in N(\Phi)$, then, by (3.34), $D(y) \subset D(\Phi)$.

Let $N^\star(\Phi)$ be the complement of $N(\Phi)$, that is,

$$W = N^\star(\Phi) \cup N(\Phi), \quad N^\star(\Phi) \cap N(\Phi) = \emptyset \tag{3.35}$$

Next, denote

$$H = S^{N(\Phi)}, \quad \Gamma_0 = S^{N^\star(\Phi)} \tag{3.36}$$

By (3.35), one can write $\Gamma = H \times \Gamma_0$ which allows us to use the general framework of the previous subsection. Let $\mathcal{B}_0$ be the $\sigma$-algebra of subsets of $\Gamma_0$ generated by the natural projections $\pi_z : \Gamma_0 \to S$, $z \in N^\star(\Phi)$. Then $H$ is a finite set and $(\Gamma_0, \mathcal{B}_0)$ is a measurable space. The projections $\pi_H : \Gamma \to H$ and $\pi_0 : \Gamma \to \Gamma_0$ and the $\sigma$-algebra $\mathcal{B}'_0$ of subsets of $\Gamma$ (which has been introduced in the paragraph a) at the beginning of Subsection 3.1) can be described in the following way:

$$\pi_H = \pi_{N(\Phi)}, \quad \pi_0 = \pi_{N^\star(\Phi)}, \quad \mathcal{B}'_0 = \mathcal{B}_{N^\star(\Phi)}, \tag{3.37}$$

returning to our original notations introduced in Section 2 in the paragraphs b) and c) just prior to the notation (2.1), while we substitute $N(\Phi)$ or $N^\star(\Phi)$ for $\Phi$. Next, using the notation (2.1), we define the following finite partitions

$$\Lambda = \Lambda_{N(\Phi)}, \quad \Lambda_1 = \Lambda_\Phi. \tag{3.38}$$

Since $\Phi \subset N(\Phi)$, it is clear that $\Lambda_1$ is the subpartition of $\Lambda$. Finally, define sub-$\sigma$-algebra

$$\mathcal{B}'_1 = \mathcal{B}_{D^\star(\Phi)}, \tag{3.39}$$

where $D^\star(\Phi)$ is the complement of $D(\Phi)$.

In summing up this paragraph, we can conclude that the notations introduced in (3.36) - (3.38) fit the general framework formulated in the previous subsection prior to Assumption 3.1. Now we can show that Assumption 3.1 itself is satisfied. More precisely, we will prove the following proposition.

**Proposition 3.5.** Let $P(x, \cdot)$ be a transition probability kernel on $(\Gamma, \mathcal{B})$, where $\Gamma = S^W$, $S$ is a finite set, $W$ is an infinite countable set, and $\mathcal{B}$ is the Borel $\sigma$-algebra of $\Gamma$, and let $\Phi$ be a given finite nonempty subset of $W$. Assume that $P(x, \cdot)$ satisfies the conditions (A1) - (A4) and define the sets $H$, $\Gamma_0$ as in (3.36). Let the projections $\pi_H$, $\pi_0$ and the $\sigma$-algebras $\mathcal{B}'_0$ and $\mathcal{B}'_1$ be specified by (3.37) and (3.39), and the partitions $\Lambda$ and $\Lambda_1$ be defined by (3.38). Then Assumption 3.1 is satisfied, that is, $\Lambda_1$ is a subpartition of $\Lambda$ and $\mathcal{B}'_1$ is a sub-$\sigma$-algebra of $\mathcal{B}'_0$, and, moreover, there exist transition probability kernels $P_H \colon H \times \Lambda_1 \to [0,1]$ and $P_0 \colon \Gamma_0 \times \mathcal{B}'_1 \to [0,1]$ such that for each $h \in H$, $B \in \Lambda_1$ we have $P_H(h, B) > 0$, and for any $x \in \Gamma$ and for each $B \in \Lambda_1$, $C \in \mathcal{B}'_1$ the formula

$$P(x, B \cap C) = P_H(\pi_H(x), B) P_0(\pi_0(x), C) \tag{3.40}$$

is satisfied.

*Proof of Proposition 3.5.* First of all, observe that, by (3.37) and (3.39) and by the fact that $D^\star(\Phi) \subset N^\star(\Phi)$, the $\sigma$-algebra $\mathcal{B}'_1$ is a sub-$\sigma$-algebra of $\mathcal{B}'_0$. Similarly, we have already pointed out (directly after the formula (3.38)) that $\Lambda_1$ is a subpartition of $\Lambda$.

Next, by (2.1) and (3.38),

$$\Lambda_1 = \{G_v \colon v \in S^\Phi\}, \tag{3.41}$$

where the sets $G_v$ have been defined in (2.2). Thus, if $B \in \Lambda_1$, one can write $B = G_v$ for some $v \in S^\Phi$. Therefore, we can define the probability kernel $P_H \colon H \times \Lambda_1 \to [0,1]$ in the following way: for each $h \in H = S^{N(\Phi)}$ and for each $B = G_v \in \Lambda_1$,

$$P_H(h, B) = \prod_{z \in \Phi} P_z\left(\pi_{N(z)}(h), \pi_z(v)\right), \tag{3.42}$$

where the probability kernels $P_z \colon S^{N(z)} \times S \to [0,1]$ have been introduced in the Assumption (A2) and the projection $\pi_{N(z)}(h)$ is well defined, since, by (3.33), $N(z) \subset N(\Phi)$.

Observe that, by the assumption **(A3)** and by (3.42), for each $h \in H$ and for each $B \in \Lambda_1$,

$$P_H(h, B) > 0.$$

Our next aim is to define the kernel $P_0 \colon \Gamma_0 \times \mathcal{B}'_1 \to [0,1]$. For that purpose, recall that, by (3.39) and by Remark 2.1, the sub-$\sigma$-algebra $\mathcal{B}'_1 = \mathcal{B}_{D^\star(\Phi)}$ is generated by the family of all the "rectangles"

$$C_w = \{u \in \Gamma \colon \pi_\Psi(u) = w\}, \tag{3.43}$$

defined for all the finite subsets $\Psi$ of $D^\star(\Phi)$ and for each $w \in S^\Psi$. Thus, it is enough to define the probability kernel $P_0$ for any such $C = C_w$, that is, for each $w \in S^\Psi$ and for any finite subset $\Psi$ of $D^\star(\Phi)$. Therefore, we can define the define the probability kernel $P_0 \colon \Gamma_0 \times \mathcal{B}'_1 \to [0,1]$ in the following way: for any $\gamma \in \Gamma_0$ and for any set $C = C_w$, where $w \in S^\Psi$ and $\Psi$ is a finite subset of $D^\star(\Phi)$,

$$P_0(\gamma, C) = \prod_{z \in \Psi} P_z\left(\pi_{N(z)}(\gamma), \pi_z(w)\right), \tag{3.44}$$

where, as in (3.42), we use the probability kernels $P_z \colon S^{N(z)} \times S \to [0,1]$ introduced in the Assumption (A2). Observe that the natural projection $\pi_{N(z)}(\gamma)$ is well defined for any given $\gamma \in \Gamma_0 = S^{N^\star(\Phi)}$, since, by Remark 3.5, if $z \in \Psi$, where $\Psi$ is a subset of $D^\star(\Phi)$, then $N(z) \cap N(\Phi) = \emptyset$, that is,

$$N(z) \subset N^\star(\Phi). \qquad (3.45)$$

(we will also use this fact later in this proof).

Also note that, by (A2), for arbitrary $x \in \Gamma$ such that $\pi_0(x) = \gamma$ one has

$$P_0(\gamma, C) = P(x, C),$$

which proves that $P_0(\gamma, \cdot)$ is really a probability kernel.

It remains to prove the formula (3.40). As we have already indicated on the previous stage, the sub-$\sigma$-algebra $\mathcal{B}'_1 = \mathcal{B}_{D^\star(\Phi)}$ is generated by the family of all the "rectangles" $C_w$ defined by (3.43) for all the finite subsets $\Psi$ of $D^\star(\Phi)$ and for each $w \in S^\Psi$. Thus, it is enough to prove (3.40) for any such $C = C_w$.

Let $\Psi$ be a finite subset of $D^\star(\Phi)$, $w \in S^\Psi$ and $C = C_w$, and let $B = G_v \in \Lambda_1$ for some $v \in S^\Phi$. Then, by (2.2) and (3.43),

$$B \cap C = \{ u \in \Gamma : \pi_\Phi(u) = v, \pi_\Psi(u) = w \} \qquad (3.46)$$

Denote $\Psi_1 = \Psi \cup \Phi$. Since $\Psi \cap \Phi = \emptyset$, there exists $g \in S^{\Psi_1}$ such that

$$\pi_\Phi(g) = v, \quad \pi_\Psi(g) = w. \qquad (3.47)$$

Thus, one can rewrite (3.46) in the form

$$B \cap C = \{ u \in \Gamma : \pi_{\Psi_1}(u) = g \} \qquad (3.48)$$

Let $x \in \Gamma$. By (3.48), substituting $\Psi_1$ instead of $\Phi$ and $g$ instead of $v$ in the formula (2.4), we obtain, using again the fact that $\Psi \cap \Phi = \emptyset$ and (3.47),

$$P(x, B \cap C) = \prod_{z \in \Psi_1} P_z\big(\pi_{N(z)}(x), \pi_z(g)\big) = \qquad (3.49)$$

$$= \prod_{z \in \Psi_1} P_z\big(\pi_{N(z)}(x), \pi_z(g)\big) = \prod_{z \in \Psi \cup \Phi} P_z\big(\pi_{N(z)}(x), \pi_z(g)\big) =$$

$$= \prod_{z \in \Phi} P_z\big(\pi_{N(z)}(x), \pi_z(g)\big) \prod_{z \in \Psi} P_z\big(\pi_{N(z)}(x), \pi_z(g)\big) =$$

$$= \prod_{z \in \Phi} P_z\big(\pi_{N(z)}(x), \pi_z(v)\big) \prod_{z \in \Psi} P_z\big(\pi_{N(z)}(x), \pi_z(w)\big)$$

However, by (3.33) and (3.37), for each $z \in \Phi$ we have: $N(z) \subset N(\Phi)$, and, therefore, for each $x \in \Gamma$,

$$\pi_{N(z)}(x) = \pi_{N(z)}\big(\pi_{N(\Phi)}(x)\big) = \pi_{N(z)}\big(\pi_H(x)\big) \qquad (3.50)$$

Hence, by (3.42) and (3.50),

$$P_H(\pi_H(x), B) = \prod_{z \in \Phi} P_z\big(\pi_{N(z)}(\pi_H(x)), \pi_z(v)\big) = \qquad (3.51)$$

$$= \prod_{z \in \Phi} P_z\big(\pi_{N(z)}(x), \pi_z(v)\big)$$

Similarly, by (3.45) and (3.37), for each $z \in \Psi$ we have: $N(z) \subset N^\star(\Phi)$, that is, for each $x \in \Gamma$,

$$\pi_{N(z)}(x) = \pi_{N(z)}\big(\pi_{N^\star(\Phi)}(x)\big) = \pi_{N(z)}(\pi_0(x)), \qquad (3.52)$$

which yields, by (3.44),

$$P_0(\pi_0(x), C) = \prod_{z \in \Psi} P_z\big(\pi_{N(z)}(\pi_0(x)), \pi_z(w)\big) = \qquad (3.53)$$

$$= \prod_{z \in \Psi} P_z\big(\pi_{N(z)}(x), \pi_z(w)\big)$$

Finally, by (3.51) and (3.53), the formula (3.49) yields (3.40).

**Q.E.D**

*Additional notations*

Recall that in the previous subsection we have introduced the $\sigma$-algebra $\mathcal{A}$ of subsets of $\Gamma \times \Gamma$ generated by the projections $\tilde{\pi}_L$ and $\tilde{\pi}_R$ with respect to $\mathcal{B}$, where $\tilde{\pi}_L(u,v) = u$ and $\tilde{\pi}_R(u,v) = v$ for any $(u,v) \in \Gamma \times \Gamma$, that is, $\mathcal{A} = \mathcal{B} \otimes \mathcal{B}$. We have also introduced the set $M(\Gamma \times \Gamma)$ of all the probability measures defined on $(\Gamma \times \Gamma, \mathcal{A})$, and the sub-$\sigma$-algebra $\mathcal{A}_0 = \mathcal{B}'_0 \otimes \mathcal{B}'_1$ of $\mathcal{A}$. At the present stage, considering the interpretations (3.37) and (3.39), it will be more convenient to write $\mathcal{A}_\Phi$ instead of $\mathcal{A}_0$, that is

$$\mathcal{A}_\Phi = \mathcal{B}'_0 \otimes \mathcal{B}'_1 = \mathcal{B}_{N^\star(\Phi)} \otimes \mathcal{B}_{D^\star(\Phi)}, \qquad (3.54)$$

(to put it another way, $\mathcal{A}_\Phi$ is the sub-$\sigma$-algebra generated by the semi-algebra of the sets of form $D \times C \in \mathcal{A}$ such that $D \in \mathcal{B}_{N^\star(\Phi)}, C \in \mathcal{B}_{D^\star(\Phi)}$.

Similarly, in the previous section we have defined the finite partition $\Delta = \Lambda \times \Lambda_1$ of $\Gamma \times \Gamma$. In the present subsection, according to the interpretation (3.38), we will write, instead of it,

$$\Delta_\Phi = \Lambda_{N(\Phi)} \times \Lambda_\Phi. \qquad (3.55)$$

Next, let $\mu \in M(\Gamma \times \Gamma)$ and let $P(x, \cdot)$ be a given transition probability kernel on $(\Gamma, \mathcal{B})$.

As in (3.3), define the measure $\mu^P \in M(\Gamma \times \Gamma)$ by the formula

$$\mu^P(A \times B) = \int_A P(x, B) \mu_L(dx) \qquad (3.56)$$

for any $A, B \in \mathcal{B}$. For the sake of brevity denote by $D_\Phi(\mu \| \mu^P)$ the relative entropy of $\mu$ with respect to $\mu^P$ corresponding to the partition $\Delta_\Phi$, that is,

$$D_\Phi(\mu \| \mu^P) = D_{\Delta_\Phi}(\mu \| \mu^P), \qquad (3.57)$$

where $D_\Delta(\mu \| \mu^P)$ has been defined by (3.8).

Now we can bring together the results of Propositions 3.4 and 3.5.

**Corollary 3.6.** Let $P(x, \cdot)$ be a transition probability kernel on $(\Gamma, \mathcal{B})$, where $\Gamma = S^W$, $S$ is a finite set, $W$ is an infinite countable set, and $\mathcal{B}$ is the Borel $\sigma$-algebra of $\Gamma$, and let $\Phi$ be a given finite nonempty subset of $W$. Assume that $P(x, \cdot)$ satisfies the conditions (A1) - (A4) and let $\mu \in M(\Gamma \times \Gamma)$ be such that

$$D(\mu \| \mu^P) < \infty, \qquad (3.58)$$

where $D(\mu \| \mu^P)$ is the relative entropy of $\mu$ with respect to $\mu^P$ ( as it has been defined by the formulae (3.4) -(3.5)). Then,

$$\sup_{G \in \mathcal{A}_\Phi} |\mu^P(G) - \mu(G)| \leq \frac{1}{\sqrt{2}} \left( D(\mu \| \mu^P) - D_\Phi(\mu \| \mu^P) \right)^{\frac{1}{2}}. \qquad (3.59)$$

where $\mathcal{A}_\Phi$ and $D_\Phi(\mu \| \mu^P)$ have been defined in (3.54) and (3.57).

*Proof.* By Proposition 3.5, the conditions of Proposition 3.4 are satisfied. By substituting (3.54) and (3.57) into (3.16), we obtain (3.59).

**Q.E.D**

Finally, we will formulate the resulting conclusion of the present subsection. For this purpose, we will need some additional notations. For any $\mu \in M(\Gamma \times \Gamma)$ we will define the left and right marginal measures $\mu_L, \mu_R \in M(\Gamma)$ by the formulae

$$\mu_L(A) = \mu(A \times \Gamma), \quad \mu_R(A) = \mu(\Gamma \times A), \qquad (3.60)$$

for any $A \in \mathcal{B}$. Introduce *the set of the measures with symmetrical marginal distributions* by

$$M_S = \{ \mu \in M(\Gamma \times \Gamma): \mu_L = \mu_R \}. \qquad (3.61)$$

Now we can formulate and prove the final result of the present subsection.

**Corollary 3.7.** Let $P(x, \cdot)$ be a transition probability kernel on $(\Gamma, \mathcal{B})$, where $\Gamma = S^W$, $S$ is a finite set, $W$ is an infinite countable set, and $\mathcal{B}$ is the Borel $\sigma$-algebra of $\Gamma$. Assume that $P(x, \cdot)$ satisfies the conditions (A1) - (A4) and let $\nu \in M(\Gamma)$ be such that

$$I(\nu) < \infty, \qquad (3.62)$$

where $I(\nu)$ is the action functional defined by (1.2). Then there exists $\mu \in M_S$ such that $\mu_L = \nu$, $D(\mu \| \mu^P) = I(\nu)$, and for any given finite nonempty subset $\Phi$ of $W$ we have

$$\sup_{A \in \mathcal{B}_{N^*(\Phi)}} |v^P(A) - v(A)| \leq \frac{1}{\sqrt{2}} \left( D(\mu \| \mu^P) - D_\Phi(\mu \| \mu^P) \right)^{\frac{1}{2}}. \tag{3.63}$$

where the measure $v^P \in M(\Gamma)$ has been introduced in (2.5), and $\mathcal{B}_{N^*(\Phi)}$ and $D_\Phi(\mu \| \mu^P)$ have been defined in (3.37) and (3.57), respectively.

*Proof.* For the convenience of the reader, let us recall some basic properties of the action functional $I(v)$ defined by (1.2). Using our notations, we can rewrite the formula (2.21) on the page 401 of the paper [5] of Donsker and Varadhan in the following way

$$I(v) = \inf_{\mu \in M_S: \mu_L = v} \bar{I}(\mu). \tag{3.64}$$

where (using again our notations), $\bar{I}(\mu)$ is given by the definition (2.4) of [5]

$$\bar{I}(\mu) = -\inf \left\{ \ln \left( \int_{\Gamma \times \Gamma} u \, d\mu^P \right) - \int_{\Gamma \times \Gamma} \ln(u) \, d\mu \right\}$$

where the infimum is taken over all positive continuous functions $u$ defined on the compact $\Gamma \times \Gamma$. However, it is a well-known fact (see for instance, Th. 5.2.1 of [11]), that for any $\mu \in M(\Gamma \times \Gamma)$ one has

$$D(\mu \| \mu^P) = \bar{I}(\mu). \tag{3.65}$$

The minimum is reached, since the set $\{\mu \in M_S: \mu_L = v\}$ is a compact and $\bar{I}(\mu)$ is a lower semi-continuous functional, which, together with (3.64) and (3.65) implies that

$$I(v) = \min_{\mu \in M_S: \mu_L = v} D(\mu \| \mu^P). \tag{3.66}$$

By (3.66), there exists $\mu \in M_S$ such that $\mu_L = \mu_R = v$, $D(\mu \| \mu^P) = I(v) < \infty$.

Therefore, all the conditions of our Corollary 3.6 are satisfied, and, thus, by (3.59),

$$\sup_{G \in \mathcal{A}_\Phi} |\mu^P(G) - \mu(G)| \leq \frac{1}{\sqrt{2}} \left( D(\mu \| \mu^P) - D_\Phi(\mu \| \mu^P) \right)^{\frac{1}{2}}, \tag{3.67}$$

where, recall, $\mathcal{A}_\Phi = \mathcal{B}_{N^*(\Phi)} \otimes \mathcal{B}_{D^*(\Phi)}$. Since for any $A \in \mathcal{B}_{N^*(\Phi)}$ and for any set of form $G = A \times \Gamma$ we have $\mu(G) = v(A)$, $\mu^P(G) = v^P(A)$, we derive immediately the estimate (3.63).

**Q.E.D**

## 3. PROOF OF THE LEMMA 2.1: THE FINAL STAGE.

Let $P(x, \cdot)$ be a transition probability kernel on $(\Gamma, \mathcal{B})$, such that the conditions (A1) - (A4) are satisfied, and let $v \in M(\Gamma)$ be such that $I(v) < \infty$. According to the conditions of the lemma, a sequence of finite subsets $\Phi_n$, $n \geq 1$ of $W$ is such that

$$\Phi_n \subset \Phi_{n+1}, \quad \bigcup_{n=1}^{\infty} \Phi_n = W. \tag{3.68}$$

Substituting $\Phi_n$ for each $n \geq 1$ instead of $\Phi$ to the formula (3.63), we can reformulate Corollary 3.7 in the following way: there exists $\mu \in M_S$ such that $\mu_L = \nu$, $D(\mu \| \mu^P) = I(\nu) < \infty$, and for any $n \geq 1$ it holds

$$\sup_{A \in \mathcal{B}_{N^\star(\Phi_n)}} |\nu^P(A) - \nu(A)| \leq \frac{1}{\sqrt{2}} \left( D(\mu \| \mu^P) - D_{\Phi_n}(\mu \| \mu^P) \right)^{\frac{1}{2}}. \qquad (3.69)$$

Recall that $\mathcal{B}_{\Phi_n^*}$ is the sub-$\sigma$-algebra of $\mathcal{B}$ generated by the projections $\pi_z : \Gamma \to S$, $z \in \Phi_n^*$, and, similarly, $\mathcal{B}_{N^\star(\Phi_n)}$ is the sub-$\sigma$-algebra of $\mathcal{B}$ generated by the projections $\pi_z : \Gamma \to S$, $z \in N^\star(\Phi_n)$, where $N^\star(\Phi_n) = W - N(\Phi_n)$ and, by (3.33),

$$N(\Phi_n) = \bigcup_{z \in \Phi_n} N(z). \qquad (3.70)$$

Denote for $n \geq 1$

$$\alpha_n = \sup_{A \in \mathcal{B}_{\Phi_n^*}} |\nu^P(A) - \nu(A)| \geq 0 \qquad (3.71)$$

By (3.68), we have $\mathcal{B}_{\Phi_{n+1}^*} \subset \mathcal{B}_{\Phi_n^*}$, and, therefore, $\alpha_n$ is a monotonically non-increasing sequence. On the other hand, by (A2), the subset $N(\Phi_n)$ is finite for any $n \geq 1$. Thus, by (3.68), for any $n \geq 1$ there exists an integer $m$ large enough such that $N(\Phi_n) \subset \Phi_m$, which yields $\mathcal{B}_{\Phi_m^*} \subset \mathcal{B}_{N^\star(\Phi_n)}$, and, therefore,

$$\alpha_m \leq \sup_{A \in \mathcal{B}_{N^\star(\Phi_n)}} |\nu^P(A) - \nu(A)| . \qquad (3.72)$$

By (3.69) and (3.72), since $\alpha_m$ is a monotonically non-increasing sequence, we have, for any $n \geq 1$,

$$0 \leq \lim_{m \to \infty} \alpha_m \leq \frac{1}{\sqrt{2}} \left( D(\mu \| \mu^P) - D_{\Phi_n}(\mu \| \mu^P) \right)^{\frac{1}{2}}. \qquad (3.73)$$

Observe that by (3.8), (3.55) and (3.57),

$$D_{\Phi_n}(\mu \| \mu^P) = \sum_{A \times B \in \Delta_{\Phi_n}} \mu(A \times B) \ln\left(\frac{\mu(A \times B)}{\mu^P(A \times B)}\right). \qquad (3.74)$$

where $\Delta_{\Phi_n} = \Lambda_{N(\Phi_n)} \times \Lambda_{\Phi_n}$. Denote by $\mathcal{A}_n$ the finite the sub-$\sigma$-algebra of $\mathcal{A} = \mathcal{B} \otimes \mathcal{B}$ generated by the partition $\Delta_{\Phi_n}, n \geq 1$, then $D_{\Phi_n}(\mu \| \mu^P)$ is the relative entropy of the measures $\mu$ and $\mu^P$ restricted to $\mathcal{A}_n$.

Next, by (3.68) and (3.70), one has $N(\Phi_n) \subset N(\Phi_{n+1})$. Thus, the partition $\Delta_{\Phi_{n+1}}$ refines the partition $\Delta_{\Phi_n}$, which yields $\mathcal{A}_n \subset \mathcal{A}_{n+1}$. On the other hand, by (3.68), the sequence of partitions $\Lambda_{\Phi_n}$, $n \geq 1$, generates the Borel $\sigma$-algebra $\mathcal{B}$ of $\Gamma$, and the same fact is true for the sequence of the partitions $\Lambda_{N(\Phi_n)}, n \geq 1$. Thus, the sequence $\mathcal{A}_n$ generates asymptotically the Borel $\sigma$-algebra $\mathcal{A}$ of $\Gamma \times \Gamma$.

We are now in the position to apply Corollary 5.2.3 of [11], which gives

$$\lim_{n \to \infty} D_{\Phi_n}(\mu \| \mu^P) = D(\mu \| \mu^P).$$

Finally, the last statement, together with (3.71) and (3.73), completes the proof.

**Q.E.D**

**DECLARATIONS**
**Competing interests**

- No funds, grants, or other support was received.

- The author has no competing interests to declare that are relevant to the content of this article.

- The author has no relevant financial or non-financial interests to disclose.

**Data Availability Statement**

No datasets were generated or analyzed during the current study.